\documentclass{amsart}
\def\S{\mathbb{S}}

\def\R{\mathbb{R}}
\def\C{\mathbb{C}}
\def\P{\mathbb{P}}

\def\Z{\mathbb{Z}}

\newtheorem{theorem}{Theorem}[section]

\newtheorem{corollary}[theorem]{Corollary}
\newtheorem{proposition}[theorem]{Proposition}

\theoremstyle{definition}

\theoremstyle{remark}
\newtheorem{remark}[theorem]{Remark}

\numberwithin{equation}{section}

\begin{document}

\title[Hamiltonian stability and index]{Hamiltonian stability and index of \\ minimal Lagrangian surfaces of \\ the complex projective plane}

\author{Francisco Urbano}
\address{Departamento de Geometr\'{\i}a y Topolog\'{\i}a, Universidad de Granada, 18071-Granada, Spain}
\email{furbano@ugr.es}
\thanks{The author was supported in part by  a MEC-Feder grant No. MTM2004-00109.}


\date{January 11, 2005}


\begin{abstract}
We show that the Clifford torus and the totally geodesic real projective plane $\R\P^2$ in the complex projective plane $\C\P^2$ are the unique Hamiltonian stable minimal Lagrangian compact surfaces of $\C\P^2$ with genus $g\leq 4$, when the surface is orientable, and with Euler characteristic $\chi\geq -1$, when the surface is nonorientable. Also we characterize $\R\P^2$ in $\C\P^2$ as the least possible index minimal Lagrangian compact nonorientable surface of $\C\P^2$.
\end{abstract}

\maketitle

\section{Introduction}
The second variation operator of minimal submanifolds of Riemannian manifolds (the {\em Jacobi operator}) carries the information about the stability properties of the submanifold when it is thought as a critical point for the volume functional. When the ambient Riemannian manifold is the complex projective space $\C\P^n$, Lawson and  Simons [LS] characterized the complex submanifolds as the unique stable minimal submanifolds of $\C\P^n$. In particular, minimal Lagrangian submanifolds of $\C\P^n$ are unstable. In [O1] Oh introduced the notion of Hamiltonian stability for minimal Lagrangian submanifolds of $\C\P^n$ (or more generally of any K\"{a}hler manifold), as those ones such that the second variation of volume is nonnegative for Hamiltonian deformations of $\C\P^n$. He proved that the Clifford torus in $\C\P^n$ is Hamiltonian stable and conjectured that it is also volume minimizing under Hamiltonian deformations of $\C\P^n$. B. Kleiner had proved that the totally geodesic Lagrangian real projective space $\R\P^n$ in $\C\P^n$ is volume minimizing under Hamiltonian deformations. 
 
In [U] Urbano (see also [H] Theorem B) characterized the Clifford torus as the unique Hamiltonian stable minimal Lagrangian torus in $\C\P^2$ and got a lower bound for the index (the number of negative eigenvalues of the Jacobi operator) of the minimal Lagrangian compact orientable surfaces in $\C\P^2$, proving that the index is at least $2$ and it is $2$ only for the Clifford torus. In this paper the author continues studying these problems, proving, among others, the following results:

\begin{quote}
{\em
The Clifford torus is the unique Hamiltonian stable minimal Lagrangian compact orientable surface of $\C\P^2$ with genus $g\leq 4$}.
\end{quote}
\begin{quote}
{\em The totally geodesic real projective plane $\R\P^2$ is the unique Hamiltonian stable minimal Lagrangian compact nonorientable surface of $\C\P^2$ with Euler characteristic $\chi\geq -1$}.
\end{quote}
\begin{quote}
{\em The index of a minimal Lagrangian compact nonorientable surface of $\C\P^2$ is at least $3$ and it is $3$ only for the totally geodesic real projective plane $\R\P^2$.}
\end{quote}
To prove these results we need to have a control of the index of the minimal Lagrangian Klein bottles of $\C\P^2$ which admit a one-parameter group of isometries. Following the ideas of [CU], in section 5 we describe explicitly those minimal surfaces and estimate their index.

The author would like to thank A. Ros for his valuable comments about the paper.
\section{Preliminaries}
Let $\C\P^2$ be the complex projective plane with its canonical Fubini-Study metric $\langle,\rangle$ of constant holomorphic sectional curvature $4$. Then 
\[
\C\P^2=\{\Pi(z)=[z]\,/\,z\in\S^5\},
\]
where $\Pi:\S^5\rightarrow\C\P^2$ is the Hopf projection, being $\S^5$ the unit sphere in the complex Euclidean space $\C^3$.
The complex structure $J$ of $\C^3$ induces via $\Pi$ the canonical complex structure $J$ on $\C\P^2$ (we will denote both by $J$). The K\"{a}hler two form in $\C\P^2$ is defined by $\omega(.,.)=\langle J.,.\rangle$. 

An immersion $\Phi:\Sigma\rightarrow\C\P^2$ of a surface $\Sigma$ is called Lagrangian if $\Phi^*\omega=0$. This means that the complex structure $J$ defines a bundle isomorphism from the tangent bundle to $\Sigma$ onto the normal bundle to $\Phi$, allowing us to identify the sections on the normal bundle $\Gamma(T^{\perp}\Sigma)$ with the $1$-forms on $\Sigma$ by
\begin{eqnarray}
\Gamma(T^{\perp}\Sigma)&\equiv& \Omega^1(\Sigma) \\
\xi&\equiv& \alpha\nonumber
 \end{eqnarray}
being $\alpha$ the 1-form on $\Sigma$ defined by $\alpha (v)=\omega(v,\xi)$ for any $v$ tangent to $\Sigma$, and where $\Omega^p(\Sigma)$, $p=0,1,2$, denotes the space of $p$-forms on $\Sigma$.

From now on our Lagrangian surface will be minimal and compact. Among them, we would like to bring out the {\em Clifford torus}
\[
T=\{[z]\in\C\P^2\,|\,|z_i|^2=\frac{1}{3},i=1,2,3\},
\]
and the {\em totally geodesic Lagrangian real projective plane}
\[
\R\P^2=\{[z]\in\C\P^2\,|\,z_i=\bar z_i,i=1,2,3\},
\]
whose $2:1$ oriented covering provides the totally geodesic Lagrangian immersion $\S^2\rightarrow\C\P^2$ of the unit sphere. An important property of these surfaces (see for instance [EGT]), which will be use in the paper , is that {\em $\R\P^2\subset\C\P^2$ is the unique minimal Lagrangian projective plane immersed in $\C\P^2$} and hence {\em $\S^2\rightarrow\C\P^2$ is the unique minimal Lagrangian sphere immersed in $\C\P^2$}. 

Using the identification (2.1), if $L:\Gamma(T^{\perp}\Sigma)\rightarrow\Gamma(T^{\perp}\Sigma)$ denotes the Jacobi operator of the second variation of the area, Oh proved (in [O1]) that $L$ is given by 
\begin{eqnarray*}
L:\Omega^1(\Sigma) \rightarrow \Omega^1(\Sigma) \\
L=\Delta +6I,
\end{eqnarray*}
where $I$ is the identity map and, in general, $\Delta=\delta d+d\delta$ is the Laplacian on $\Sigma$ acting on $p$-forms, $p=0,1,2$, being $\delta$ the codifferential operator of the exterior differential $d$. Hence, {\it the index of $\Phi$, that we will denote by $\hbox{Ind}\,(\Sigma)$, is the number of eigenvalues (counted with multiplicity) of $\Delta:\Omega^1(\Sigma) \rightarrow \Omega^1(\Sigma)$ less than $6$}.

To study the Jacobi operator, we  consider the Hodge decomposition
\[
\Omega^1(\Sigma)={\mathcal H}(\Sigma)\oplus d\,\Omega^0(\Sigma)\oplus\delta\Omega^2(\Sigma),
\]  
which allows to write, in a unique way, any 1-form $\alpha$ as $\alpha=\alpha_0 + df +\delta \beta$, being $\alpha_0$ a harmonic 1-form, $f$ a real function and $\beta$ a 2-form on $\Sigma$. The space of harmonic 1-forms, ${\mathcal H}(\Sigma)$, is the kernel of $\Delta$ and its dimension is the first Betti number of $\Sigma$: $\beta_1(\Sigma)$. As $\Delta$ commutes with $d$ and $\delta$, the positive eigenvalues of $\Delta:\Omega^1(\Sigma) \rightarrow \Omega^1(\Sigma)$ are the positive eigenvalues of $\Delta:\Omega^0(\Sigma)\rightarrow\Omega^0(\Sigma)$ joint to the positive eigenvalues of $\Delta:\Omega^2(\Sigma)\rightarrow\Omega^2(\Sigma)$. Hence 
\begin{equation}
\hbox{Ind}\,(\Sigma)=\beta_1(\Sigma)\,+\,\hbox{Ind}_0(\Sigma)\,+\,\hbox{Ind}_1(\Sigma),
\end{equation}
where {\it $\hbox{Ind}_0(\Sigma)$ is the number of positive eigenvalues (counted with multiplicity) of $\Delta:\Omega^0(\Sigma)\rightarrow \Omega^0(\Sigma)$  less than $6$} and {\it $\hbox{Ind}_1(\Sigma)$ is the number of positive eigenvalues (counted with multiplicity) of $\Delta:\Omega^2(\Sigma)\rightarrow \Omega^2(\Sigma)$ less than $6$}.  When $\Sigma$ is $\R\P^2\subset\C\P^2$, $\hbox{Ind}_0(\Sigma)=0$, $\hbox{Ind}_1(\Sigma)=3$ and hence $\hbox{Ind}(\Sigma)=3$.

If the compact surface $\Sigma$ is orientable, the star operator $\star:\Omega^0(\Sigma)\rightarrow\Omega^2(\Sigma)$ says us that the eigenvalues of $\Delta$ acting on $\Omega^0(\Sigma)$ or on $\Omega^2(\Sigma)$ are the same, and so $\hbox{Ind}_0(\Sigma)=\hbox{Ind}_1(\Sigma)$. Hence if $\Sigma$ is a minimal Lagrangian compact {\it orientable} surface of $\C\P^2$ of genus $g$, then 
\begin{equation}
\hbox{Ind}\,(\Sigma)=2\,g\,+\,2\,\hbox{Ind}_0(\Sigma).
\end{equation}
When $\Sigma$ is the totally geodesic Lagrangian two-sphere in $\C\P^2$, $\hbox{Ind}_0(\Sigma)=3$ and $\hbox{Ind}(\Sigma)=6$.

The variational vector fields of the Hamiltonian deformations of the Lagrangian surface $\Sigma$ are the normal components of the Hamiltonian vector fields on $\C\P^2$. If $X=J\bar\nabla F$, for certain smooth function $F:\C\P^2\rightarrow\R$, is a Hamiltonian vector field on $\C\P^2$, the $1$-form associated to the normal component of $X$, under the identification (2.1), is $d(F\circ\Phi)$. So our minimal Lagrangian compact surface $\Sigma$ is Hamiltonian stable if the first positive eigenvalue of $\Delta$ acting on $\Omega^0(\Sigma)$ is at least $6$. But it is well-known that $6$ is always an eigenvalue of $\Delta:\Omega^0(\Sigma)\rightarrow\Omega^0(\Sigma)$ (see proof of Theorem 3.3). Hence {\it $\Sigma$ is Hamiltonian stable if the first positive eigenvalue of $\Delta$ acting on $\Omega^0(\Sigma)$ is $6$}. In this setting it is natural to call to $\hbox{Ind}_0(\Sigma)$ the {\em Hamiltonian index} of $\Sigma$.

Let $\Sigma$ be a nonorientable Riemannian compact surface and $\pi:\tilde{\Sigma}\rightarrow\Sigma$ the $2:1$ orientable Riemannian covering. If $\tau:\tilde{\Sigma}\rightarrow\tilde{\Sigma}$ is the change of sheet involution, then the spaces of forms on $\tilde{\Sigma}$ can be decompose in the following way:
\[
\Omega^i(\tilde{\Sigma})=\Omega^i_+(\tilde{\Sigma})\oplus\Omega^i_-(\tilde{\Sigma}),\quad i=0,1,2,
\]
where
\[
\Omega^i_{\pm}(\tilde{\Sigma})=\{\alpha\in\Omega^i(\tilde{\Sigma})\,/\,\tau^*\alpha=\pm \alpha\}.
\]
Also the space of harmonic $1$-forms on $\tilde{\Sigma}$ decomposes into two subspaces ${\mathcal H}(\tilde{\Sigma})={\mathcal H}_+(\tilde{\Sigma})\oplus {\mathcal H}_-(\tilde{\Sigma})$, where again ${\mathcal H}_{\pm}(\tilde{\Sigma})=\{\alpha\in{\mathcal H}(\tilde{\Sigma})\,/\,\tau^*\alpha=\pm\alpha\}$. In this way we obtain
\[
\Omega^1_{\pm}(\tilde{\Sigma})={\mathcal H}_{\pm}(\tilde{\Sigma})\oplus d\,\Omega^0_{\pm}(\tilde{\Sigma})\oplus \delta\,\Omega^2_{\pm}(\tilde{\Sigma}).
\]
As $\pi\circ\tau=\pi$, the map $\alpha\in\Omega^i(\Sigma)\mapsto\pi^*\alpha\in\Omega^i(\tilde{\Sigma})$ allows to identify ${\mathcal H}(\Sigma)\equiv{\mathcal H}_+(\tilde{\Sigma})$ and $\Omega^i(\Sigma)\equiv\Omega^i_+(\tilde{\Sigma})$, $i=0,1,2$. Also, as $\Sigma$ is nonorientable, $\star\tau^*=-\tau^*\star$, and so $\star$ identifies $\Omega^0_-(\tilde{\Sigma})\equiv\Omega^2_+(\tilde{\Sigma})$. Hence we obtain the identification 
\begin{eqnarray*}
\Omega^2(\Sigma)&\equiv&\Omega^0_-(\tilde{\Sigma})\\
\alpha&\equiv& f  
\end{eqnarray*}
where $\pi^*\alpha=f\omega_0$, being $\omega_0$ the volume $2$-form on $\tilde{\Sigma}$.

Now, let $\Phi:\Sigma\rightarrow\C\P^2$ be a minimal Lagrangian immersion of a compact nonorientable surface $\Sigma$ and $\Phi\circ\pi:\tilde{\Sigma}\rightarrow\C\P^2$ the corresponding minimal Lagrangian immersion of its 2:1 orientable covering $\tilde{\Sigma}$. As $\Sigma$ is nonorientable, the eigenvalues of $\Delta:\Omega^2(\Sigma)\rightarrow\Omega^2(\Sigma)$ are positives, and hence, taking into account the above remarks, {\em $\hbox{Ind}_1(\Sigma)$ is the number of eigenvalues (counted with multiplicity) of $\Delta:\Omega^0_-(\tilde{\Sigma})\rightarrow\Omega^0_-(\tilde{\Sigma})$ less than $6$}. Also, as {\em $\hbox{Ind}_0(\Sigma)$ is the number of positive eigenvalues (counted with multiplicity) of $\Delta:\Omega^0_+(\tilde{\Sigma})\rightarrow\Omega^0_+(\tilde{\Sigma})$ less than $6$}, we obtain that
\begin{equation}
\hbox{Ind}_0(\Sigma)+\hbox{Ind}_1(\Sigma)=\hbox{Ind}_0(\tilde{\Sigma)},
\end{equation}
and hence from (2.3)
\[
2\hbox{Ind}(\Sigma)=\hbox{Ind}(\tilde{\Sigma}).
\]

\section{Proof of the results}
\begin{theorem}
Let $\Phi:\Sigma\rightarrow\C\P^2$ be a Hamiltonian stable minimal Lagrangian immersion of a compact orientable surface of genus $g$. If $g\leq 4$ then $\Phi$ is an embedding and $\Phi(\Sigma)$ is the Clifford torus.
\end{theorem}
{\it Proof:\/} As $\Sigma$ is Hamiltonian stable, the first positive eigenvalue of $\Delta:\Omega^0(\Sigma)\rightarrow\Omega^0(\Sigma)$ is $6$. Now we use a well-known argument. From the Brill-Noether theory, we can get a nonconstant meromorphic map $\phi:\Sigma\rightarrow\S^2$ of degree $d\leq 1+[\frac{g+1}{2}]$, where $[.]$ stands for integer part. Then there exists a Moebius transformation $F:\S^2\rightarrow\S^2$ such that $\int_{\Sigma}(F\circ\phi)=0$, and so
\[
\int_{\Sigma}|\nabla(F\circ\phi)|^2\geq 6\int_{\Sigma}|F\circ\phi|^2=6\,\hbox{Area}(\Sigma).
\]
But $\int_{\Sigma}|\nabla(F\circ\phi)|^2=8\pi\,\hbox{degree}(F\circ\phi)=8\pi\,\hbox{degree}(\phi)\leq 8\pi(1+[\frac{g+1}{2}])$. Hence we obtain that $3\hbox{Area}(\Sigma)\leq 4\pi(1+[\frac{g+1}{2}])$.

On the other hand, Montiel and Urbano ([MU] Corollary 6) proved that $\hbox{Area}(\Sigma)\geq 2\pi\mu$ ($\mu$ being the maximum multiplicity of the immersion $\Phi$) and that the equality holds if and only if the surface is the totally geodesic two-sphere. So we obtain that 
\[
3\mu\leq 2(1+[\frac{g+1}{2}]), 
\]
and the equality implies that the surface is the totally geodesic Lagrangian two-sphere.

Using [EGT], Lemma 4.6, we know that $\Phi$ is not an embedding (i.e. $\mu\geq 2$) when $g\geq 2$. So in this case $2<[\frac{g+1}{2}]$, which is a contradiction when $g=2,3,4$. If $g=0$, our surface is the totally geodesic Lagrangian two-sphere, which is Hamiltonian unstable. Hence the surface must be a torus and using [U], Corollary 2, we conclude that it is the Clifford torus.\hfill q.e.d.

\begin{remark}
If $\lambda_1$ is the first positive eigenvalue of $\Delta:\Omega^0(\Sigma)\rightarrow\Omega^0(\Sigma)$ and $g\geq2$ then the above reasoning proves that $\lambda_1<2(1+[\frac{g+1}{2}])$. Hence if $g=2$ we obtain that $\lambda_1<4$.
\end{remark}

\begin{theorem}
Let $\Phi:\Sigma\rightarrow\C\P^2$ be a minimal Lagrangian immersion of a compact nonorientable surface $\Sigma$ with Euler characteristic $\chi(\Sigma)\geq -1$. Then 
\begin{enumerate}
\item If $\Sigma$ is a projective plane with a handle ($\chi(\Sigma)=-1$) then $\Sigma$ is Hamiltonian unstable.
\item If $\Sigma$ is a Klein bottle ($\chi(\Sigma)=0$) then $\hbox{Ind}_0(\Sigma)\geq 2$. In particular $\Sigma$ is  Hamiltonian unstable.
\end{enumerate}
As consequence, if $\Phi$ is Hamiltonian stable then $\Phi$ is an embbeding and $\Phi(\Sigma)$ is $\R\P^2$.
\end{theorem}
{\it Proof:\/} 
 We will denote also by $\langle,\rangle$ the Euclidean metric in $\C^3$. In the Lie algebra $\mathrm s\mathrm o(6)$ of the isometry group of $\S^5$, we consider the subspace
$\mathrm s\mathrm o^+(6)=\{A\in \mathrm s\mathrm o(6)\,/\,AJ=JA\quad\hbox{and}\quad \hbox{Trace}\,AJ=0\}$, which is the real representation of the Lie algebra ${\mathrm s\mathrm u}(3)$. Then for any $A\in \mathrm s\mathrm o^+(6)$, the function on the sphere $p\in\S^5\mapsto\langle Ap,Jp\rangle\in\R$ can be projected to $\C\P^2$, defining a map
\begin{eqnarray*}
F_A:\C\P^2\rightarrow\R\\
F_A(\Pi(p))=\langle Ap,Jp\rangle
\end{eqnarray*}
First we compute the gradient of $F_A$. If $v$ is any vector tangent to $\C\P^2$ at $\Pi(p)$, then 
\[
v\cdot F_a=2\langle Av^*,Jp\rangle,
\]
being $v^*$ the horizontal lifting of $v$ to $T_p\S^5$. So
\[
(\bar\nabla F_A)_{\Pi(p)}=-2(d\Pi)_p(AJp+F_A(p)p),
\]
for any $\Pi(p)\in\C\P^2$. Taking derivatives again and using that $\Pi:\S^5\rightarrow\C\P^2$ is a Riemannian  submersion, one has that the Hessian of $F_A$ is given by
\begin{equation}
(\bar\nabla^2F_A)(v,w)=-2F_A\langle v,w\rangle+2\langle Av^*,Jw^*\rangle,
\end{equation}
for any vectors $v,w\in T_{\Pi(p)}\C\P^2$. 

Now, if $\Phi:\Sigma\rightarrow\C\P^2$ is a minimal Lagrangian immersion of a compact surface $\Sigma$ and $f_A:\Sigma\rightarrow\R$ is defined by $f_A=F_A\circ\Phi$, then by decomposition 
\[
\bar\nabla F_A=\nabla f_A +\xi
\]
in its tangential and normal components and taking into account (3.1) we deduce
\[
(\nabla^2f_A)(v,w)=-2f_A\langle v,w\rangle+2\langle Av^*,Jw^*\rangle+\langle\sigma(v,w),\xi\rangle,
\]
which implies that $\Delta f_A+6f_A=0$. So we have defined a linear map
\begin{eqnarray*}
H:\mathrm s\mathrm o^+(6)&\rightarrow& V_6=\{f\,/\,\Delta f+6f=0\}\\
A&\mapsto& f_A,
\end{eqnarray*}
and hence the multiplicity of the eigenvalue $6$ satisfies $m(6)\geq 8-\hbox{dim\,Ker}\,H$. If $A\in\hbox{Ker}H$, then $f_A=0$ and so $\nabla f_A=0$, which implies that $\bar\nabla F_A=\xi$. Using (3.1), the tangent vector field $J\xi$ satisfies 
\[
\langle\nabla_v J\xi,v\rangle=2\langle Av^*,v^*\rangle=0,
\]
which means that $J\xi$ is a Killing field on $\Sigma$. If $J\xi=0$, then $\bar{\nabla}F_A$ vanishes identically on the points of the surface, which implies, looking at the expression of $\bar{\nabla}F_A$, that $A=0$. Hence $\hbox{dim\,Ker}\,H\leq \hbox{dim\,\{Killing fields on}\, \Sigma\}$. Finally we get that
\begin{equation}
m(6)\geq8-\hbox{dim\,\{Killing fields on}\, \Sigma\}.
\end{equation}
In that follows we will use the following Nadirashvili's result

\vspace{0.1cm}

{\sc Theorem A} [N] {\em Let $\Sigma$ be a compact nonorientable surface with Euler characteristic $\chi(\Sigma)\leq 0$ and $\langle,\rangle$ any Riemannian metric on $\Sigma$. Then the multiplicity of the $i$-th eigenvalue of the Laplacian $\Delta:\Omega^0(\Sigma)\rightarrow\Omega^0(\Sigma)$ satisfies $m(\lambda_i)\leq 3+2i-2\chi(\Sigma)$.}

\vspace{0.1cm}

First suppose that $\Sigma$ is a Hamiltonian stable projective plane with a handle, i.e. $\chi(\Sigma)=-1$. Then the first positive eigenvalue of $\Delta:\Omega^0(\Sigma)\rightarrow\Omega^0(\Sigma)$ is $6$. But Theorem A says that $m(6)=m(\lambda_1)\leq 7$. So (3.2) implies that there exists a non-trivial Killing vector field in our surface, which is impossible. This proves part 1.

Suppose now that $\Sigma$ is a Hamiltonian stable Klein bottle, i.e. $\chi(\Sigma)=0$. Again the first positive eigenvalue of $\Sigma$ is $\lambda_1=6$ and Theorem A says that $m(6)=m(\lambda_1)\leq 5$. From (3.2), $\hbox{dim\,\{Killing fields on}\, \Sigma\}\geq 3$, which is impossible. 

Suppose now that $\Sigma$ is a Klein bottle with $\hbox{Ind}_0(\Sigma)=1$. Then $\lambda_1<6$, the multiplicity of $\lambda_1$ is $m(\lambda_1)=1$ and $\lambda_2=6$. Using again Theorem A, $m(6)=m(\lambda_2)\leq 7$. From (3.2), our Klein bottle admits a nontrivial Killing field. Proposition 5.1 (see section 5) says that $\Sigma$ is congruent to some finite Riemannian covering of $K_{n,m}$ for certain integers $n,m$. Then, from Proposition 5.2, we have that
\[
\hbox{Ind}_0(\Sigma)\geq\hbox{Ind}_0(K_{n,m})\geq 6,
\] 
which is a contradiction. This proves part 2.

Finally, if $\Sigma$ is Hamiltonian stable, then $\Sigma$ is a projective plane, i.e. $\chi(\Sigma)=1$, and hence our surface is $\R\P^2$. This finishes the proof. 
\hfill
q.e.d.

\begin{theorem}
Let $\Phi:\Sigma\rightarrow\C\P^2$ be a minimal Lagrangian immersion of a Klein bottle or a projective plane with a handle. Then $\hbox{Ind}_1(\Sigma)\geq 1$.
\end{theorem}
{\it Proof:\/} Let $\pi:\tilde{\Sigma}\rightarrow\Sigma$ be the $2:1$ orientable Riemannian covering of $\Sigma$ and $\tau:\tilde{\Sigma}\rightarrow\tilde{\Sigma}$ the change of sheet involution. If $\hbox{Ind}_1(\Sigma)=0$ then, taking into account the last remarks of section 2, the first eigenvalue $\lambda_1$ of $\Delta:\Omega^0_-(\tilde{\Sigma})\rightarrow\Omega^0_-(\tilde{\Sigma})$ satisfies $\lambda_1\geq 6$. Hence 
\[
\int_{\tilde{\Sigma}}|\nabla f|^2\geq 6\int_{\tilde{\Sigma}}f^2,\quad \forall f\in C^{\infty}(\tilde{\Sigma})\quad\hbox{such that}\quad f\circ\tau=-f.
\]
From Theorem 1 in [RS], we can get a nonconstant meromorphic map $\phi:\tilde{\Sigma}\rightarrow\S^2$ satisfying $\phi\circ\tau=-\phi$ of degree $d\leq 1+g$, where $g$ is the genus of the compact orientable surface $\tilde{\Sigma}$. Hence we obtain 
\[
\int_{\tilde{\Sigma}}|\nabla \phi|^2\geq 6\int_{\tilde{\Sigma}}|\phi|^2=6\,\hbox{Area}(\tilde{\Sigma}).
\]
But $\int_{\tilde{\Sigma}}|\nabla \phi|^2=8\pi\,\hbox{degree}(\phi)\leq 8\pi(1+g)$. So we get that $3\hbox{Area}(\tilde{\Sigma})\leq 4\pi(1+g)$. Now, as $\hbox{Area}(\tilde{\Sigma})=2\hbox{Area}(\Sigma)$, using again Corollary 6 in [MU] as in the proof of Theorem 3.1, we have that $3\mu\leq 1+g$ ($\mu$ being the maximum multiplicity of $\Phi$) and the equality implies that $\Sigma$ is $\R\P^2$.

As our surface is either a Klein bottle or a projective plane with a handle, we have that $g=1\,\hbox{or}\,2$ and that the equality (in the above inequality) is not attained, i.e. $3\mu<1+g$. This is a contradiction and the proof is finished.                                                                                                \hfill q.e.d.
\begin{remark}
In this nonorientable case, we cannot use that the maximum multiplicity $\mu$ of $\Phi$ satisfies $\mu\geq 2$ when the Euler characteristic $\chi(\Sigma)\leq 0$. The author only knows that $\Phi$ is not an embedding when $\Sigma$ is a Klein bottle, i.e. $\chi(\Sigma)=0$. (See [M], Theorem 2).
\end{remark}
\begin{corollary}
Let $\Phi:\Sigma\rightarrow\C\P^2$ be minimal Lagrangian immersion of a compact nonorientable surface $\Sigma$. Then $\hbox{Ind}(\Sigma)\geq3$ and the equality holds if and only if $\Phi$ is an embedding and $\Phi(\Sigma)$ is $\R\P^2$.
\end{corollary}
{\it Proof:\/} If $\Sigma$ is a projective plane, then as we mentioned in section 2, $\Phi$ is an embedding and $\Phi(\Sigma)$ is $\R\P^2$, whose index is $3$.

If $\Sigma$ is a compact nonorientable surface with Euler characteristic $\chi(\Sigma)\leq 0$, then $\beta_1(\Sigma)=1-\chi(\Sigma)$. So, from (2.2) $\hbox{Ind}(\Sigma)\geq\beta_1(\Sigma)\geq4$ when $\chi(\Sigma)\leq -3$. If $\Sigma$ is either a Klein bottle or a projective plane with a handle (i.e. $\chi(\Sigma)=0$ or $-1$), Theorems 3.3 and 3.4 joint with (2.2) say again that $\hbox{Ind}(\Sigma)\geq 4$. Finally, if $\Sigma$ is a projective plane with $2$ handles, i.e. $\chi(\Sigma)=-2$, and $\tilde{\Sigma}$ is its $2:1$ orientable covering, then the genus of $\tilde{\Sigma}$ is $3$, and  Theorem 3.1 joint with (2.4) say that 
\[
\hbox{Ind}(\Sigma)=3+\hbox{Ind}_0(\Sigma)+\hbox{Ind}_1(\Sigma)=3+\hbox{Ind}_0(\tilde{\Sigma})\geq 4.
\]
This finishes the proof.\hfill
q.e.d.
\section{Area minimizing surfaces in their Hamiltonian isotopy classes.}
A Lagrangian immersion $\Phi:\Sigma\rightarrow\C\P^2$ is called {\em Hamiltonian minimal} [O2] if it is a critical point of the area functional for Hamiltonian deformations. The corresponding Euler-Lagrange equation says that $\Sigma$ is Hamiltonian minimal if and only if $\hbox{div}\,JH=0$, where $H$ is the mean curvature vector of $\Sigma$ and $\hbox{div}$ stands for the divergence operator on $\Sigma$. 

Suppose now that $\Phi:\Sigma\rightarrow\C\P^2$ is a minimal Lagrangian immersion of a compact surface $\Sigma$ and that in its Hamiltonian isotopy class there exists a minimizer $\tilde{\Sigma}$ for the area. Then, in particular, $\tilde{\Sigma}$ is a Hamiltonian minimal Lagrangian surface, and so $\hbox{div}\,J\tilde H=0$, being $\tilde{H}$ the mean curvature vector of $\tilde{\Sigma}$. But, using Theorem I in [O3], the deRham cohomology class defined by the mean curvature vector $\alpha(v)=\omega(v,\tilde H)$ is invariant under Hamiltonian isotopies. As $\Sigma$ is minimal, this class must vanish, and hence there exists a smooth function $f:\tilde{\Sigma}\rightarrow\R$ such that $\nabla f=J\tilde H$. Now, $\hbox{div}\,J\tilde H=0$ implies that $f$ is a harmonic function and so it is constant, which means that $\tilde{\Sigma}$ is also minimal. Moreover, as $\tilde{\Sigma}$ is also Hamiltonian stable, our previous results say that:
\begin{quote}
{\em The Clifford torus of $\C\P^2$ is the unique area minimizing surface in its Hamiltonian isotopy class, provided that someone existed.} 
\end{quote}
\begin{quote}
{\em There exists no area minimizing surfaces in the Hamiltonian isotopy class of a minimal Lagrangian compact orientable surface of $\C\P^2$ of genus $2$,$3$ or $4$.}
\end{quote}
\begin{quote}
{\em There exists no area minimizing surfaces in the Hamiltonian isotopy class of a minimal Lagrangian Klein bottle or a minimal Lagrangian projective plane with a handle of $\C\P^2$.}
\end{quote}

\section{Minimal Lagrangian Klein bottles of $\C\P^2$ admitting a one-parameter group of isometries.}

In this section we are going to describe the minimal Lagrangian Klein bottles of $\C\P^2$ admitting a one-parameter group of isometries, estimating also their Hamiltonian index. To understand it, we need to give a short introduction to the elliptic Jacobi functions. We will follow the notation and the results of [CU]. 

Given $p\in [0,1[$, let $\hbox{dn}(x,p)=\hbox{dn}(x)$, $\hbox{cn}(x,p)=\hbox{cn}(x)$ and $\hbox{sn}(x,p)=\hbox{sn}(x)$ be the elementary Jacobi elliptic functions with modulus $p$. Then, the following properties are well known:
\begin{equation}
\hbox{sn}^2(x)+\hbox{cn}^2(x)=1,\quad \hbox{dn}^2(x)+p^2\hbox{sn}^2(x)=1, \quad\forall x\in\R
\end{equation}
and
\begin{eqnarray}
\frac{d}{dx}\hbox{dn}(x)&=&-p^2\hbox{sn}(x)\hbox{cn}(x),\nonumber\\
\frac{d}{dx}\hbox{cn}(x)&=&-\hbox{sn}(x)\hbox{dn}(x),\\
\frac{d}{dx}\hbox{sn}(x)&=&\hbox{cn}(x)\hbox{dn}(x)\nonumber.
\end{eqnarray} 
Also, if 
\[
K(p)=\int_0^{\pi/2}\frac{d\theta}{\sqrt{1-p^2\sin^2\theta}}
\]
is the {\em complete elliptic integral of the first kind}, then the elliptic functions have the following symmetry and periodicity properties:
\begin{eqnarray}
\hbox{dn}(x+2K)&=&\hbox{dn}(x),\quad\hbox{dn}(K-x)=\hbox{dn}(K+x),\nonumber\\
\hbox{cn}(x+2K)&=&-\hbox{cn}(x),\quad\hbox{cn}(K-x)=-\hbox{cn}(K+x),\\
\hbox{sn}(x+2K)&=&-\hbox{sn}(x),\quad\hbox{sn}(K-x)=\hbox{sn}(K+x)\nonumber.
\end{eqnarray}
In particular all of them are periodic of period $4K$ and $\hbox{dn},\,\hbox{cn}$ are even, i.e. $\hbox{dn}(-x)=\hbox{dn}(x),\,\hbox{cn}(-x)=\hbox{cn}(x)$, meanwhile $\hbox{sn}$ is odd, i.e. $\hbox{sn}(-x)=-\hbox{sn}(x)$. Moreover $\hbox{dn}$ is a positive function, $\hbox{cn}(x)$ vanishes for $x=(2k+1)K,\,k\in\Z$ and $\hbox{sn}(x)$ vanishes for $x=2kK,\,k\in\Z$. 

Although in [CU] the authors classified the (non-totally geodesic) minimal Lagrangian immersions in $\C\P^2$ of simply-connected surfaces invariants by a one-parameter group of isometries of $\C\P^2$, in fact they only used that condition in order to have a non-trivial Killing field on the surface. So, really, they classified the minimal Lagrangian immersions in $\C\P^2$ of simply-connected surfaces admitting a Killing field.

Let $\Phi:(\R^2,g)\rightarrow\C\P^2$ be a minimal Lagrangian isometric immersion such that $(\R^2,g)$ admits a Killing field. Then (see [CU]) the metric $g$ can be written as $g=e^{2u(x)}(dx^2+dy^2)$, where $u$ is a solution of the following problem
\begin{equation}
u''(x)+e^{2u(x)}-e^{-4u(x)}=0,\quad e^{2u(0)}=b\in[1,\infty[,\quad u'(0)=0.
\end{equation}
We will denote the metric by $g_b$. The solutions of (5.4) are given by
\begin{equation}
e^{2u(x)}=b(1-q^2\hbox{sn}^2(rx,p)),
\end{equation}
where
\[
q^2=1-\frac{1+\sqrt{1+8b^3}}{4b^3},\, r^2=b-\frac{1-\sqrt{1+8b^3}}{4b^2},\, p^2=\frac{bq^2}{r^2}.
\]
Hence the solutions $u(x)$ of (5.4) are periodic functions with period $2K/r$  and satisfy $u(-x)=u(x), \forall x\in\R$. The only constant solution of (5.4) corresponds to $b=1$ and the associated minimal Lagrangian immersion is the universal covering of the Clifford torus. 

On the other hand, in [CU], Theorem 4.1, the minimal Lagrangian immersions of $(\R^2,g_b)$ into $\C\P^2$ were explicitly given. Using a reasoning like in [CU], Theorem 4.2, it is not difficult to prove that the minimal Lagrangian immersion $(\R^2,g_b)\rightarrow\C\P^2$ corresponding to the initial condition $b$ is the universal covering of a minimal Lagrangian Klein bottle in $\C\P^2$ if and only if the number $\frac{1+\sqrt{1+8b^3}}{4b^3}$, which belongs to the interval $]0,1]$, satisfies the two following conditions:
\begin{enumerate}
\item $\frac{1+\sqrt{1+8b^3}}{4b^3}$ is a rational number, and 
\item if $\frac{1+\sqrt{1+8b^3}}{4b^3}=\frac{m}{n}$ with $m,n\in\Z$, and $(m,n)=1$, then $n$ is odd and $m+2n=\dot6$.
\end{enumerate}
We note that $0<m<n$ and that $(2n+m)/3$ is an even integer, and $(n+2m)/3$ and $(n-m)/3$ are odd integers.   
Moreover, in such case, the corresponding group $G_{n,m}$ of transformations of $\R^2$ is generated by 
\[
(x,y)\mapsto(x+4K/r,y),\quad (x,y)\mapsto(-x,y+\sqrt 2mb\pi/3)
\]
 and the corresponding minimal Lagrangian immersion is given by
\begin{eqnarray*}
\Psi_{n,m}&:&\R^2\rightarrow\C\P^2\\
\Psi_{n,m}(x,y)&=&\left[\left(\lambda\,\hbox{dn}(rx)e^{\frac{i(n+m)y}{\sqrt2mb}},\mu\,\hbox{cn}(rx)e^{\frac{-iny}{\sqrt2mb}},\nu\,\hbox{sn}(rx)e^{\frac{-iy}{\sqrt2b}}\right)\right],
\end{eqnarray*}
where 
\[
\lambda^2=\frac{n}{2n+m},\quad \mu^2=\frac{n+m}{2n+m},\quad \nu^2=\frac{n+m}{n+2m},\quad r^2=\frac{n(n+2m)}{2b^2m^2},
\]
and where the modulus of the elliptic Jacobi functions is given by $p^2=(n^2-m^2)/n(n+2m)$.

If $K_{n,m}=\R^2/G_{n,m}$ is the associated Klein bottle and $P:\R^2\rightarrow K_{n,m}$ the projection, then the induced immersion
\begin{eqnarray*}
\Phi_{n,m}:K_{n,m}&\rightarrow&\C\P^2\\
            P(x,y)&\mapsto& \Psi_{n,m}(x,y)
\end{eqnarray*}
defines a minimal Lagrangian immersion of the Klein bottle $K_{n,m}$ in $\C\P^2$.

We can summarize the above reasoning in the following result.
\begin{proposition}
Let $\Phi:\Sigma\rightarrow\C\P^2$ be a minimal Lagrangian immersion of a Klein bottle $\Sigma$ admitting a one-parameter group of isometries. Then $\Phi$ is congruent to some finite Riemannian covering of $\Phi_{n,m}:K_{n,m}\rightarrow\C\P^2$ with $n$ and $m$ integers such that
$0<m<n$, $(m,n)=1$, $n$ is odd and $2n+m=\dot 6$, where $\dot 6$ stands for the positive integer multiples of $6$. 
\end{proposition}
Now we use a similar method to use by Haskins in [H], Theorem E, in order to estimate the Hamiltonian index of the Klein bottles $K_{n,m}$.
Following the proof of Theorem 3.3, the eigenspace of $\Delta:\Omega^0(K_{n,m})\rightarrow\Omega^0(K_{n,m})$ corresponding to the eigenvalue $6$ has at least dimension $7$, and the functions $\{g_i:K_{n,m}\rightarrow\R,\,1\leq i\leq 7\}$ defined by
\[
\begin{array}{ll}
g_1(P(x,y))= e^{2u(x)}-(1+2b^3)/3b^2,&   \\
g_2(P(x,y))=(\hbox{dn}\cdot \hbox{cn})(rx)\cos(\frac{2n+m}{\sqrt 2mb}y);& g_3(P(x,y))=(\hbox{dn}\cdot \hbox{cn})(rx)\sin(\frac{2n+m}{\sqrt 2mb}y),\\
g_4(P(x,y))=(\hbox{dn}\cdot \hbox{sn})(rx)\cos(\frac{n+2m}{\sqrt 2mb}y);& g_5(P(x,y))=(\hbox{dn}\cdot \hbox{sn})(rx)\sin(\frac{n+2m}{\sqrt 2mb}y),\\
g_6(P(x,y))=(\hbox{cn}\cdot \hbox{sn})(rx)\cos(\frac{n-m}{\sqrt 2mb}y);& g_7(P(x,y))=(\hbox{cn}\cdot \hbox{sn})(rx)\sin(\frac{n-m}{\sqrt 2mb}y),
\end{array}
\]
are a basis of such $7$-dimensional subspace.

Now suppose that $6=\lambda_j$ for some positive integer $j$. Then the Courant nodal Theorem (see [Ch]) says that the number of nodal sets $n_i$ of the eigenfunction $g_i$, $1\leq i\leq 7$, satisfies $n_i\leq j+1$. So to estimate $j$ we are going to compute the number of nodal sets of $g_i, 1\leq i\leq 7$. To do that we will determine the set of zeroes of $f_i=g_i\circ P:D\rightarrow\R,\,1\leq i\leq 7$ on the fundamental domain 
\[
D=\{(x,y)\in\R^2\,|\,0\leq x\leq 4K/r,\,0\leq y\leq \sqrt 2mb\pi/3\}
\]
of the Klein bottle $K_{n,m}$.

It is clear that there exists $a\in ]0,K/r[$ such that                                                            
\[
f_1^{-1}(0)=\{(x,y)\in D\,|\,x=a,\frac{2K}{r}-a,\frac{2K}{r}+a,\frac{4K}{r}-a;\,1\leq y\leq \sqrt 2mb\pi/3\}.
\]
Hence the number of nodal sets of $g_1$ is $n_1=3$.

  The zeroes of the function $(\hbox{dn}\hbox{cn})(rx)$ in the interval $[0,4K/r]$ are $K/r$ and $3K/r$. So
\begin{eqnarray*}
f_2^{-1}(0)=\{(x,y)\in D\,|\,x=K/r,3K/r;\,0\leq y\leq \sqrt 2mb\pi/3\}\cup\\
\{(x,y)\in D\,|\,0\leq x\leq 4K/r;y=\frac{\sqrt 2mb\pi}{2n+m}k,\,k=0,1,\dots,(2n+m)/3\}.
\end{eqnarray*}
Hence the number of nodal sets of $g_2$ is $n_2=2(2n+m)/3$. A similar reasoning proves that the numbers of nodal sets of the function $g_3$ is $n_3=2(2n+3)/3$. 

In a similar way, the zeroes of the function $(\hbox{dn}\hbox{sn})(rx)$ in the interval $[0,4K/r]$ are $0$, $2K/r$ and $4K/r$. So
\begin{eqnarray*}
f_4^{-1}(0)=\{(x,y)\in D\,|\,x=0,2K/r,4K/r;\,0\leq y\leq \sqrt 2mb\pi/3\}\cup\\
\{(x,y)\in D\,|\,0\leq x\leq 4K/r;y=\frac{\sqrt 2mb\pi}{n+2m}k,\,k=0,1,\dots,(n+2m)/3\}.
\end{eqnarray*}
Hence the number of nodal sets of $g_4$ (and of $g_5$) is $n_4=n_5=2(n+2m)/3$.

Finally, the zeroes of the function $(\hbox{cn}\hbox{sn})(rx)$ in the interval $[0,4K/r]$ are $0,K/r$, $2K/r$, $3K/r$ and $4K/r$. So
\begin{eqnarray*}
f_6^{-1}(0)=\{(x,y)\in D\,|\,x=0,K/r,2K/r,3K/r,4K/r;\,0\leq y\leq \sqrt 2mb\pi/3\}\cup\\
\{(x,y)\in D\,|\,0\leq x\leq 4K/r;y=\frac{\sqrt 2mb\pi}{n-m}k,\,k=0,1,\dots,(n-m)/3\}.
\end{eqnarray*}
Hence the number of nodal sets of $g_6$ (and of $g_7$) is $n_6=n_7=4(n-m)/3$.

Hence, as $2n+m>n+2m$ and $2n+m>2(n-m)$ we have obtained that if $6=\lambda_j$ then 
\[
j+1\geq 2(2n+m)/3\geq 8.
\]
Using that fact joint with Theorem 3.4 and (2.2) we finally obtain the following result.
\begin{proposition}
Let $\Phi_{n,m}:K_{n,m}\rightarrow\C\P^2$ be the minimal Lagrangian immersion of the Klein bottle $K_{n,m}$ with $n,m$ integers satisfying $0<m<n$, $(n,m)=1$, $n$ odd
and $2n+m=\dot 6$. Then
\[
\hbox{Ind}_0(K_{n,m})\geq \frac{2(2n+m)}{3}-2\geq 6,\quad \hbox{Ind}(K_{n,m})\geq \frac{2(2n+m)}{3}\geq 8.
\]
\end{proposition}
To finalize, we are going to study other interesting properties of the Klein bottles family $\{K_{n,m}\}$.
\begin{proposition}
Let $\Phi_{n,m}:K_{n,m}\rightarrow\C\P^2$ be the minimal Lagrangian immersion of the Klein bottle $K_{n,m}$ with $n,m$ integers satisfying $0<m<n$, $(n,m)=1$, $n$ odd
and $2n+m=\dot 6$. Then
\begin{enumerate}
\item The first eigenvalue of $\Delta:\Omega^0(K_{n,m})\rightarrow\Omega^0(K_{n,m})$ satisfies
\[
\lambda_1(K_{n,m})<2-\frac{1}{2b^3}=2-\frac{m^2}{n(n+m)},
\]
\item The area of the Klein bottle $K_{n,m}$ is given by
\[
A(K_{n,m})=\frac{4\pi\sqrt n}{3\sqrt{n+2m}}((n+2m)E-mK),
\]
\end{enumerate}
where $E=\int_0^{\pi/2}\sqrt{1-p^2\sin^2\theta}\,d\theta$ is the {\em complete elliptic integral of the second kind}.  
\end{proposition}
{\it Proof:\/}
Let $f:K_{n,m}\rightarrow\R$ be the function defined by $f(P(x,y))=\hbox{sn}(rx)$. Then from (5.3) and (5.5) we have that
\begin{eqnarray*}
\int_{K_{n,m}}f\,dA=\int_{0}^{4K/r}\int_0^{\sqrt2mb\pi/3}\hbox{sn}(rx)e^{2u(x)}dydx\\=\frac{\sqrt2mb^2\pi}{3}\int_0^{4K/r}(\hbox{sn}(rx)-q^2\hbox{sn}^3(rx))dx=0.
\end{eqnarray*}
Hence, if $\lambda_1$ is the first eigenvalue of $\Delta:\Omega^0(K_{n,m})\rightarrow\Omega^0(K_{n,m})$, 
\[
-\int_{K_{n,m}}f\Delta fdA\geq\lambda_1\int_{K_{n,m}}f^2dA.
\]
But using (5.1),(5.2),(5.3) and (5.5) we have   
\[
\Delta f(P(x,y))=e^{-2u(x)}\frac{d^2}{dx^2}\,\hbox{sn}(rx)=-2f(P(x,y))+\frac{e^{-2u(x)}}{2b^2}f(P(x,y)),
\]
and so, using that $e^{-2u(x)}\geq 1/b$, we have that
\[
\lambda_1\int_{K_{n,m}}f^2dA<2\int_{K_{n,m}}f^2dA-\frac{1}{2b^3}\int_{K_{n,m}}f^2dA,
\]
which proves (1).

On the other hand, from (5.1) and (5.5) it follows that
\[
A(K_{n,m})=\int_{0}^{4K/r}\int_0^{\sqrt2mb\pi/3}e^{2u(x)}dydx=\frac{\sqrt2mb\pi}{3}\int_{0}^{4K/r}(b-r^2+r^2\hbox{dn}^2(rx))dx.
\]
If $E(u)=\int_{0}^u\hbox{dn}^2(y)\,dy$, then it is known that $E(u+2K)=E(u)+2E$, $\forall u\in \R$. Using this property in the above expression we obtain (2). 

\hfill q.e.d.




\end{document}